\documentclass[twoside,bezier,reqno]{amsart}

\usepackage{epsfig}
\usepackage{amsmath,amsthm,amsfonts,amssymb,amscd,euscript}

\input{epsf}
\newcommand{\filebegin}{\begin{document}}
\newcommand{\fileend}{\end{document}}
\makeatother
\def\thefootnote{}
\newcommand{\lo}{\longrightarrow}
\newcommand{\NMM}{\hspace*{2mm}}

\renewcommand{\baselinestretch}{1.1}

\input{epsf}
\renewcommand{\baselinestretch}{1.1}
\def\n{\noindent}%

\numberwithin{equation}{section}
\def\mapdown#1{\Big\downarrow\rlap
{$\vcenter{\hbox{$\scriptstyle#1$}}$}}
\newtheorem{theorem}{Theorem}[section]
\newtheorem{lemma}[theorem]{Lemma}
\newtheorem{proposition}[theorem]{Proposition}
\newtheorem{corollary}[theorem]{Corollary}

\theoremstyle{definition}
\newtheorem{definition}[theorem]{Definition}
\newtheorem{example}[theorem]{\sc Example}
\newtheorem{xca}[theorem]{Exercise}
\theoremstyle{remark}
\newtheorem{remark}[theorem]{Remark}
\begin{document}
%%%%%%%%%%%%%%%%%%%%%%%%%%%%%%%%%%%%%%%

\setcounter{page}{1} \noindent

%%%%%%%%%%%%%%%%%%%%%%%%%%%%%%%%%%%%%%%
\vspace*{2cm}
\begin{center}
{\bf\large  On the  Gorensetein $(n,d)$-Flat and Gorensetein $(n,d)$-Injective Modules}
 \\[0.5cm]
{M. Amini$^{^{*}}$ \footnote{$^*$Corresponding Author} \\[2mm]
Department of Mathematics, Payame Noor University, Tehran, Iran\\

{\tt E-mail: mostafa.amini@pnu.ac.ir}\\

}
\end{center}
\vspace*{0.5cm}

%------------------------------------------------------------------------------------------------------------------------------
\begin{quote}
{\small \hfill{\rule{13.3cm}{.1mm}\hskip2cm} \textbf{Abstract.}
 Let $R$ be a ring. In this paper,  Gorensetein $(n,d)$-flat and Gorenstein $(n,d)$-injective modules are introduced and some of their basic properties are studied. Moreover, some characterizations of rings over   Gorensetein $(n,d)$-flat and Gorenstein $(n,d)$-injective modules are given. Also, some known results can be obtained as corollaries.  }\\[.2cm]
%\vspace{1mm} {\renewcommand{\baselinestretch}{1}
%\parskip = 10 mm
{
\noindent{\small {\it \bf 2010 MSC}\,:13D07; 16D40; 18G25;}}\\
\noindent{\small {\it \bf Keywords}\,:Gorenstein $(n,d)$-Flat module, Gorensetein $(n,d)$-Injective module, $n$-Coherent ring.}\\
\vspace{-3mm}\hfill{\rule{13.3cm}{.1mm}\hskip2cm}
\end{quote}
%_

\markboth 
{M. Amini}
 {On the  Gorensetein $(n,d)$-Flat and Gorensetein $(n,d)$-Injective Modules}

%%%%%%%%%%%%%%%%%%%%%%%%%%%%%%%%%%%%%%%%%%%%%%%%%%%%%%%%%%%%%%%%%%%%%%%%%%%%%%%%%%%%%%%%%%%%%%%%%%%%%%%%%%%%%%%%%%
%%%%%%%%%%%%%%%%%%%%%%%%%%%%%%%%%%%%%%%%%%%%%%%%%%%%%%%%%%%%%%%%%%%%%%%%%%
\section{Basic Definitions and Notations}
In this paper, we have assumed that $R$ is an  associative ring with non-zero identity and  modules are
unitary. In this section, first some fundamental concepts and notations are stated.
 Let  $n,d$ non-negative integers and $M$ a right $R$-module. Then 
\begin{enumerate}
\item[(1)]
$M$ is said to be $n$-{\it presented}  \cite{X.J} if there is an exact sequence of  right $R$-modules
 $ F_{n}\rightarrow F_{n-1}\rightarrow\dots\rightarrow F_1\rightarrow F_0\rightarrow
M\rightarrow$ , where each $F_i$ is a finitely generated free.
\item[(2)]
$R$ is called right $n$-{\it coherent} ring \cite{ X.J} if every  $n$-presented  right $R$-module
is $(n + 1)$-presented.
\item[(3)] 
$M$ is called  $(n,d)$-flat \cite{ N.I, X.J} if ${Tor}_{d+1}^{R}(U,M)=0$ for every $n$-presented left $R$-module $U$.
\item[(4)] 
 $M$ is called  $(n,d)$-injective  \cite{ N.I, X.J} if ${ Ext}_{R}^{d+1}(U,M)=0$ for every $n$-presented right $R$-module $U$.
\item[(5)]
 $M$  is said to be {\it Gorenstein flat} (resp.,{\it Gorenstein injective})
\cite{E.J,E.B,L.X} if there
is an exact sequence $\cdots\rightarrow I_1\rightarrow I_{0}\rightarrow I^0\rightarrow
I^1\rightarrow\cdots$ of flat (resp., injective ) right modules
with $M= {\rm ker}(I^0\rightarrow
I^1)$ such that $U\otimes_{R}-$ (resp; ${ Hom}(U,-)$) leaves the sequence exact whenever
$U$ is an  left injective (resp; right injective)  module.
\item [(6)] For any commutative ring $R$ of prime characteristic $p > 0$, assume that $F_ R : R \rightarrow R^ {(e)}$ is the $e$-th iterated Frobenius map in which $R^{(e)}\cong R$. 
Then, the {\it perfect closure} \cite{k.h} of $R$, denoted  by $R_{\infty}$, is defined as the limit of the following direct system:
\[
\begin{CD}
R @>F_R>> R @>F_R>> R @>F_R>> \cdots \\
\end{CD}
\]
\item [(7)] Assume that $S \geq R$ is a unitary ring extension. Then, 
the ring $S$ is called  {\it $R$-projective} \cite{X.J,wu} in case, for any  right $S$-module $M_S$ with
an right $S$-module $N_S$, $N_R\mid M_R$ implies $N_S\mid M_S$, where $N\mid M$ means that $N$ is a direct summand
of $M$.
%\item [(8)]  The ring extension $S \geq R$ is called a {\it finite normalizing extension} \cite{X.J,wu} in case there is a finite subset $\{s_1,\cdots,s_n\}\subseteq S$ such that
%$S=\sum_{i=1}^{i=n}s_{i}R$ and $s_{i}R=Rs_{i}$ for $i=1,\cdots,n$.
\item [(8)] A finite normalizing extension $S \geq R$ is called an {\it almost excellent extension} \cite{X.J,wu} in
case $_RS$ is flat, $S_R$ is projective, and the ring $S$ is $R$-projective. 
\end{enumerate}

For more information about the relative cohence of rings and modules, see \cite{Z.G,cam,k.h,aejm1}.
 In this paper, we introduce the {Gorenstein $(n,d)$-flat} modules and { Gorenstein $(n,d)$-injective}
 modules. A right module $G$ is said to be {\it Gorenstein $(n,d)$-flat} ({\it $G_{n}^{d}$-flat}
for short), if there exists the following exact sequence of $(n,d)$-flat right $R$-modules
$${\mathbf{F}}= \cdots\longrightarrow F_1\longrightarrow F_{0}\longrightarrow F^0\longrightarrow
F^1\longrightarrow\cdots$$ with $G=\ker(F^0\rightarrow F^1)$ such that ${\rm Tor}_{d}^{R}(U,{\mathbf{F}})$ leaves this sequence exact whenever $U$ is an $n$-presented left $R$-module with ${fd}_{R}(U)<\infty$.
A right module $G$ is said to be {\it Gorenstein $(n,d)$-injective} ({\it $G_n^{d}$-injective}
for short), if there exists the following exact sequence of $(n,d)$-injective right $R$-modules
$${\mathbf{A}}= \cdots\longrightarrow A_1\longrightarrow A_{0}\longrightarrow A^0\longrightarrow
A^1\longrightarrow\cdots$$ with $G=\ker(A^0\rightarrow A^1)$ such that ${\rm Ext}_{R}^{d}(U,{\mathbf{A}})$ leaves this sequence exact whenever $U$ is an $n$-presented right $R$-module with ${pd}_{R}(U)<\infty$.
The Gorenstein $(n,d)$-flat dimension (resp; Gorenstein $(n,d)$-injective dimension) of
a module $M$ is denoted by $G_{n}^{d}.fd(M)$ (resp; $G_{n}^{d}.id(M)$). $G_{n}^{d}.fd(M)$ of an $R$-module $M$
is defined that $G_{n}^{d}.fd(M)\leq m$ (resp; $G_{n}^{d}.id(M)\leq m$) if and only if $M$ has a Gorenstein $(n,d)$-flat (resp; Gorenstein $(n,d)$-injective)
resolution of length $m$.

In Section 2, we study some basic properties of the Gorenstein $(n,d)$-flat and Gorenstein $(n,d)$-injective modules.
In section 3, we give  sufficient conditions under which every module is Gorenstein $(n,d)$-injective. For instance, if $R$ is a $n$-coherent ring, then every $R$-module is Gorenstein $(n,d)$-injective if and only if $R$ is $(n,d)$-injective  if and only if every Gorenstein $(n,0)$-injective is Gorenstein $(n,d)$-flat.  Finally in section 4, some results of  Gorenstein (n, d)-flat  and
Gorenstein (n, d)-injective modules on projective algebras
 are given.  For instance, if $f: R\rightarrow S$ be a surjective ring homomorphism,  $S$ a
projective $R$-module and $G$ a  right $S$-module, then $G_R$ is Gorenstein $(n,d)$-flat (resp; $G_{R}$ is Gorenstein (n, d)-injective )  if and only if $G_S$ is Gorenstein $(n,d)$-flat (resp; $G_{S}$ is Gorenstein (n, d)-injective ) if and only if $S\otimes_{R}G$ (resp; ${\rm Hom}_{R}(S,G)$) is Gorenstein $(n,d)$-flat (resp; Gorenstein $(n,d)$-injective  ). Also, if $S\geq R$ is an almost excellent extension, then the { perfect closure}  of $R$ is right Gorenstein $(n,d)$-injective if and only if the { perfect closure}  of $S$ is right Gorenstein $(n,d)$-injective.

%%%%%%%%%%%%%%%%%%%%%%%%%%%%%%%%%%%%%%%%%%%%%%%%%%%%%%%%%%%%%%%%%%%%%%%%%%%%%%%%%%%%%%%%%%%%%%%%%%%%%%%%
%%%%%%%%%%%%%%%%%%%%%%%%%%%%%%%%%%%%%%%%%%%%%%%%%%%%%%%%%%%%%%%%%%%%%%%%%%%%%%%%%%%%%%%%%%%%%%%%%%%%%%%%%%%%%%%%%%%%%%%%%%%%%%%%%%%%%%%%%%
\section{Gorensetein $(n,d)$-Flatnes and Gorensetein $(n,d)$-Injectivity}
We start with the following definition.
\begin{definition}\label{2.76}
Let $G$ be a right $R$-module. Then
\begin{enumerate}
\item [\rm (1)]
 $G$  is called  Gorenstein $(n,d)$-flat if there exists the following exact sequence of $(n,d)$-flat right $R$-modules:
$${\mathbf{F}}= \cdots\longrightarrow F_1\longrightarrow F_{0}\longrightarrow F^0\longrightarrow
F^1\longrightarrow\cdots$$ with $G=\ker(F^0\rightarrow F^1)$ such that ${\rm Tor}_{d}^{R}(U,{\mathbf{F}})$ leaves this sequence exact whenever $U$ is $n$-presented left $R$-module with ${\rm fd}_R(U)<\infty$.
\item [\rm (2)]
 $G$  is called  Gorenstein $(n,d)$-injective if there exists the following  exact sequence of  $(n,d)$-injective right $R$-modules:
$${\mathbf{A}}= \cdots\longrightarrow A_1\longrightarrow A_{0}\longrightarrow A^0\longrightarrow
A^1\longrightarrow\cdots$$ with $G=\ker(A^0\rightarrow A^1)$ such that ${\rm Ext}_{R}^{d}(U,{\mathbf{A}})$ leaves this sequence exact whenever $U$  is $n$-presented right $R$-module with ${\rm pd}_R(U)<\infty.$
\end{enumerate}
\end{definition}
 
\begin{lemma}\label{2.09}
If $R$ is a left $n$-coherent ring and $U$ is an $n$-presented left $R$-module  with ${\rm fd}_{R}(U)<\infty$, then ${\rm pd}_{R}(U)<\infty$.
\end{lemma}
\begin{proof}
{ Let ${\rm fd}_{R}(U)=m<\infty$, then we must
show that ${\rm pd}_{R}(U)\leq m$. Since $R$ is a  left $n$-coherent ring and $U$ is an $n$-presented left $R$-module,  the projective resolution 
$\cdots \rightarrow F_{n+1}\rightarrow F_{n}\rightarrow\cdots \rightarrow F_{0}\rightarrow U\rightarrow 0,$ where any $F_i$ is finitely generated free, exists. On the other hand, above exact squence  is a flat resolution.  So by \cite[Proposition 8.17]{rotman},  $(m-1)$-syzygy is flat. Hence, the exact sequence
$0\rightarrow K_{m-1}\rightarrow F_{m-1}\rightarrow\cdots \rightarrow F_{0}\rightarrow U\rightarrow 0$ is a flat resolution. If $n\geq m$ or $n\leq m$,  then $K_{m-1}$ is finitely presented and consequencly by \cite[Theorem 3.56]{rotman}, $K_{m-1}$ is projective and so, ${\rm pd}_{R}(U)\leq m$. }
\end{proof}

In the following theorem, we show that in the case of $n$-coherent rings, the existence of  complexs  of a module is sufficient to be Gorenstein $(n,d)$-flat and Gorenstein $(n,d)$-injective.
\begin{theorem}\label{2.3}
Let $R$ be a ring and G a  right $R$-module. 
\begin{enumerate}
\item [\rm (1)]
If $R$ is a left  $n$-coherent ring, then $G$ is Gorenstein $(n,d)$-flat if and only if there is an exact sequence 
$${\mathbf{F}}= \cdots\longrightarrow F_1\longrightarrow F_{0}\longrightarrow F^0\longrightarrow
F^1\longrightarrow\cdots$$ of  $(n,d)$-flat right $R$-modules  such that $G=\ker(F^0\rightarrow F^1)$.
\item [\rm (2)]

If $R$ is a right $n$-coherent ring, then $G$ is Gorenstein $(n,d)$-injective if and only if there is an exact sequence
$${\mathbf{A}}= \cdots\longrightarrow A_1\longrightarrow A_{0}\longrightarrow A^0\longrightarrow
A^1\longrightarrow\cdots$$
of $(n,d)$-injective right $R$-modules such that $G=\ker(A^0\rightarrow A^1)$.
\end{enumerate}
\end{theorem}
\begin{proof}
{(1) ($\Longrightarrow$) : This is a direct consequence of the definition.

($\Longleftarrow$) :  By definition, it suffices to show that ${\rm Tor}_{d}^{R}(U,{\mathbf{F}})$ is exact for
every  $n$-presented left $R$-module  $U$  with ${\rm fd}_{R} (U)< \infty$. By Lemma \ref{2.09}, ${\rm pd}_{R} (U)< \infty$. We use the induction on $d$. Let $d=0$ and ${\rm pd}_{R} (U) = m < \infty$, then we  show that $U\otimes_{R} {\mathbf{F}}$ is exact.
 To prove this, we use the induction on $m$. The case $m =0$ is clear. Assume that $m \geq 1$. There exists an exact sequence
$0\rightarrow K\rightarrow P_{0}\rightarrow U \rightarrow 0$, where $P_{0}$ is projective. Now, from the $n$-coherence of $R$ and \cite[Theorem 1]{wu}, we deduce that $K$ is left $n$-presented. Also, ${\rm pd}_{R}(K) \leq m-1$. So, the following short exact sequence of complexes exists:
\begin{center}
$
\begin{array}{ccccccccc}
& \vdots&\vdots &\vdots&\\
& \downarrow & \downarrow &\downarrow & \\
0 \longrightarrow &K\otimes_{R}F_1 & \longrightarrow P_{0}\otimes_{R}F_1&\longrightarrow U\otimes_{R}F_1\longrightarrow 0 \\
& \downarrow & \downarrow &\downarrow & \\
0 \longrightarrow &K\otimes_{R}F_0 & \longrightarrow P_{0}\otimes_{R}F_0&\longrightarrow U\otimes_{R}F_0\longrightarrow 0 \\
& \downarrow &\downarrow &\downarrow & \\
0 \longrightarrow &K\otimes_{R}F^0 & \longrightarrow P_{0}\otimes_{R}F^0&\longrightarrow U\otimes_{R}F^0\longrightarrow 0 \\
& \downarrow &\downarrow &\downarrow & \\
0 \longrightarrow &K\otimes_{R}F^1 & \longrightarrow P_{0}\otimes_{R}F^1&\longrightarrow U\otimes_{R}F^1\longrightarrow 0 \\
& \downarrow &\downarrow &\downarrow & \\
& \vdots&\vdots &\vdots&\\
& \parallel & \parallel &\parallel& \\
0 \longrightarrow &K\otimes_{R}{\mathbf{F}}& \longrightarrow P_{0}\otimes_{R}{\mathbf{F}}&\longrightarrow U\otimes_{R}{\mathbf{F}}\longrightarrow 0. \\

\end{array}
$
\end{center}
\noindent By induction, $P_{0}\otimes_{R}{\mathbf{F}}$ and $K\otimes_{R}{\mathbf{F}}$ are exact, hence $U\otimes_{R}{\mathbf{F}}$ is exact by \cite[Theorem 6.10]{rotman}. 

Let $d \geq 1$. Since, $U$ is $n$-presented, the exact sequence $0\rightarrow K^{'}\rightarrow P_{0}^{'}\rightarrow U \rightarrow 0$ with $P_{0}^{'}$ is projective exists. Therefore, the following short exact sequence of complexes exists:
\begin{center}
$
\begin{array}{ccccccccc}
\vdots &\vdots&\\
\downarrow &\downarrow & \\
0 \longrightarrow  {\rm Tor}_{d}^{R}(U,F_1)&\longrightarrow {\rm Tor}_{d-1}^{R}(K^{'},F_1)\longrightarrow 0 \\
 \downarrow &\downarrow & \\
0 \longrightarrow {\rm Tor}_{d}^{R}(U,F_0)&\longrightarrow {\rm Tor}_{d-1}^{R}(K^{'},F_0)\longrightarrow 0 \\\downarrow &\downarrow & \\
0 \longrightarrow {\rm Tor}_{d}^{R}(U,F^0)&\longrightarrow {\rm Tor}_{d-1}^{R}(K^{'},F^0)\longrightarrow 0 \\
\downarrow &\downarrow & \\
0 \longrightarrow {\rm Tor}_{d}^{R}(U,F^1)&\longrightarrow {\rm Tor}_{d-1}^{R}(K^{'},F^1)\longrightarrow 0 \\
\downarrow &\downarrow & \\
\vdots &\vdots&\\
 \parallel &\parallel& \\
0 \longrightarrow {\rm Tor}_{d}^{R}(U,{\mathbf{F}})&\longrightarrow {\rm Tor}_{d-1}^{R}(K^{'},{\mathbf{F}})\longrightarrow 0. \\

\end{array}
$
\end{center}

\noindent We know that $K^{'}$ is $n$-presented with ${\rm pd}_{R} (K^{'})< \infty$, and so by induction, ${\rm Tor}_{d-1}^{R}(K^{'},{\mathbf{F}})$ is exact. Therefore, ${\rm Tor}_{d}^{R}(U,{\mathbf{F}})$ is exact and hence, $G$ is Gorenstein $(n,d)$-flat.

(2) ($\Longrightarrow$) : This is a direct consequence of the definition.

($\Longleftarrow$) Let $U$ be a  right $n$-presented $R$-module  with ${\rm pd}_{R} (U)< \infty$. Then, a similar proof to that of (1), ${\rm Ext}_{R}^{d}(U,{\mathbf{A}})$ is exact and hence $G$ is Gorenstein $(n,d)$-injective.
}

\end{proof}

\begin{remark}\label{2}
Let $R$ be a ring. Then:
\begin{enumerate}
\item [\rm (1)]
Every flat right $R$-module is $(n,d)$-flat.
\item [\rm (2)]
Every injective right $R$-module is $(n,d)$-injective.
\item [\rm (3)]
 Every right $R$-module is  flat if and only if $(1,0)$-flat.
\item [\rm (4)]
Every right $R$-module is  injective if and only if  $(0,0)$-injective.
\item [\rm (5)]
Every $(m,d)$-injective right $R$-module is $(n,d)$-injective, for any $n\geq m.$
\item [\rm (6)]
Every $(m,d)$-flat right $R$-module is $(n,d)$-flat, for any $n\geq m.$
\item [\rm (7)]
Every $n$-presented left (resp; right) $R$-module is $m$-presented, for any $n\geq m.$
\item [\rm (8)]
Every $(n,d)$-flat right $R$-module is  Gorenstein $(n,d)$-flat.
\item [\rm (9)]
Every $(n,d)$-injective right $R$-module is  Gorenstein $(n,d)$-injective.
\end{enumerate}
\end{remark}
\begin{corollary}\label{2.ty}
Let $R$ be a left $n$-coherent ring and G a right $R$-module. Then the
following assertions are equivalent:
\begin{enumerate}
\item [\rm (1)]
$G$ is Gorenstein $(n,d)$-flat;
\item [\rm (2)]
There is an exact sequence $0\rightarrow G\rightarrow B^{0}\rightarrow B^{1}\rightarrow \cdots$ of right $R$-modules, where every $B^i$ is $(n,d)$-flat;
\item [\rm (3)]
There is a short exact sequence $0\rightarrow G\rightarrow M\rightarrow L\rightarrow 0$ of right $R$-modules, where
$M$ is $(n,d)$-flat and $L $ is Gorenstein $(n,d)$-flat.
\end{enumerate}
\end{corollary}
\begin{proof}
{
$(1)\Longrightarrow (2)$ and $(1)\Longrightarrow (3)$ follow from definition.

$(2)\Longrightarrow (1)$ For $R$-module $G$, there is an exact sequence
$$\cdots\longrightarrow P_{1}\longrightarrow P_{0}\longrightarrow G\longrightarrow 0, $$
where any $P_{i}$ is flat  for any $i\geq 0$.  By Remark \ref{2}, every $P_i$ is $(n,d)$-flat. Thus, the exact sequence
$$ \cdots\longrightarrow P_1\longrightarrow P_{0}\longrightarrow B^0\longrightarrow
B^1\longrightarrow\cdots$$
of  $(n,d)$-flat right modules  exists, where $G={\rm ker}(B^0\rightarrow B^1)$. Therefore by Theorem \ref{2.3}, $G$ is Gorenstein $(n,d)$-flat, 

$(3)\Longrightarrow (2)$ Assume that the exact sequence
$$ 0\longrightarrow G\longrightarrow M \longrightarrow L\longrightarrow 0 \ \ (1)$$ exists, where
$M$  is $(n,d)$-flat and $L $ is Gorenstein $(n,d)$-flat. Since $L$ is Gorenstein $(n,d)$-flat, there is an exact sequence
$$0\longrightarrow L\longrightarrow (F^0)^{'}\longrightarrow (F^1)^{'}\longrightarrow (F^2)^{'}\longrightarrow \cdots\ \ (2)$$
where every $(F^i)^{'}$ is $(n,d)$-flat.
Assembling the sequences $(1)$ and $(2)$, we
get the exact sequence
$$0\rightarrow G\rightarrow M\rightarrow (F^0)^{'}\rightarrow (F^1)^{'}\rightarrow (F^2)^{'}\rightarrow \cdots,$$ where $M$ and any $ (F^i)^{'}$ are $(n,d)$-flat , as desired.
}
\end{proof}
\begin{corollary}\label{2.h}
Let $R$ be a right $n$-coherent ring and G a right $R$-module. Then the
following assertions are equivalent:
\begin{enumerate}
\item [\rm (1)]
$G$ is Gorenstein $(n,d)$-injective;
\item [\rm (2)]
There is an exact sequence $\cdots\rightarrow A_1\rightarrow A_{0}\rightarrow G\rightarrow 0$ of right $R$-modules, where every $A_i$ is $(n,d)$-injective;
\item [\rm (3)]
There is a short exact sequence $0\rightarrow L\rightarrow M\rightarrow G\rightarrow 0$ of $R$-modules, where
$M$ is $(n,d)$-injective and $L $ is Gorenstein $(n,d)$-injective.
\end{enumerate}
\end{corollary}
\begin{proof}
{
$(1)\Longrightarrow (2)$ and $(1)\Longrightarrow (3)$ follow from definition.

$(2)\Longrightarrow (1)$ For any module $G$, there is an exact sequence
$$0\longrightarrow G\longrightarrow I^{0}\longrightarrow I^1\longrightarrow \cdots$$
where every $I^{i}$ is injective for any $i\geq 0$. By Remark \ref{2}, each $I^{i}$ is $(n,d)$-injective.  So, the exact sequence
$$ \cdots\longrightarrow A_1\longrightarrow A_{0}\longrightarrow I^0\longrightarrow
I^1\longrightarrow\cdots$$
of $(n,d)$-injective right modules  exists, where $G={\rm ker}(I^0\rightarrow I^1)$. Therefore, $G$ is Gorenstein $(n,d)$-injective, by Theorem \ref{2.3}.

$(3)\Longrightarrow (2)$ Assume that the exact sequence
$$ 0\longrightarrow L\longrightarrow M \longrightarrow G\longrightarrow 0 \ \ (1)$$ exists, where
$M$ is $(n,d)$-injective and $L $ is Gorenstein $(n,d)$-injective. Since $L$ is Gorenstein $(n,d)$-injective, there is an exact sequence
$$\cdots\longrightarrow A_2^{'}\longrightarrow A_1^{'}\longrightarrow A_0^{'}\longrightarrow L\longrightarrow 0\ \ (2)$$
where every $A_i^{'}$ is $(n,d)$-injective.
Assembling the sequences $(1)$ and $(2)$, we
get the exact sequence
$$\cdots\rightarrow A_2^{'}\rightarrow A_1^{'}\rightarrow A_0^{'}\rightarrow M\rightarrow G\rightarrow 0,$$ where $M$ and $A_i^{'}$ are $(n,d)$-injective,
 as desired.}
\end{proof}

\begin{proposition}\label{2.5}
%For any module $G$, the following statements hold.
%\begin{enumerate}
%\item [\rm (1)]
%If $G$ is Gorenstein $(n,d)$-flat, then ${\rm Tor}_{i}^{R}(U, G)=0$
%for any $i>d$ and every $n$-presented $R$-module $U$ with ${\rm f.dim}(U) <\infty$.
%\item [\rm (2)]
%If $0\rightarrow N\rightarrow G_{m-1}\rightarrow G_{m-2}\rightarrow\cdots \rightarrow G_{0}\rightarrow G\rightarrow 0$ is an exact sequence of
%modules where every $G_i$ is a Gorenstein $(n,d)$-flat, then ${\rm Tor }_{i}^{R}(U, N)={\rm Tor}_{m+i}^{R}(U, G)$ for any $i>d$ with ${\rm f.dim} (U) <\infty$.
%\item [\rm (3)]
If $G$ is Gorenstein $(n,d)$-injective right $R$-module, then ${\rm Ext}_{R}^{i}(U, G)=0$
for any $i>d$ and every $n$-presented right $R$-module $U$  with ${\rm pd}_{R}(U) <\infty$.
%\item [\rm (4)]
%If $0\rightarrow G\rightarrow G_{0}\rightarrow G_{1}\rightarrow\cdots \rightarrow G_{m-1}\rightarrow N\rightarrow 0$ is an exact sequence of
%modules where every $G_i$ is a Gorenstein $(n,d)$-injective, then ${\rm Ext}_{R}^{i}(U, N)={\rm Ext}_{R}^{m+i}(U, G)$  for any $i>d$ with ${\rm p.dim} (U) <\infty$.

%\end{enumerate}
\end{proposition}
\begin{proof}
{
%(1) Let $G$ be a Gorenstein $(n,d)$-flat $R$-module, and ${\rm f.dim}(U)=m<\infty$.
%Then by hypothesis, the following $(n,d)$-flat resolution of $G$ exists:
%$$0\rightarrow G\rightarrow F^{0}\rightarrow \cdots \rightarrow F^{m-1}\rightarrow N\rightarrow 0.$$
%So, ${\rm Tor}_{i}^{R}(U, F^j)= 0$ for every $0 \leq j \leq m-1$ and any $i>d$, since $U$ is $n$-presented and any $F^j$ is $(n,d)$-flat. Thus, we deduce that
 %${\rm Tor}_{i}^{R}(U,G)\cong{\rm Tor}_{m+i}^{R}(U, N)=0$. 

%(2) Setting $G_m=N$ and $K_j=\ker (G_{j}\rightarrow G_{j-1})$, for every $0\leq j\leq m$, the short exact sequence
%$0\rightarrow K_{j}\rightarrow G_{j}\rightarrow K_{j-1}\rightarrow 0$ exist. Thus by (1), the induced exact sequences%
%$$0={\rm Tor}_{i+1}^{R}(U,G_{j})\rightarrow{\rm Tor}_{i+1}^{R}(U,K_{j-1})\rightarrow{\rm Tor}_{i}^{R}(U,K_{j})\rightarrow
%{\rm Tor}_{i}^{R}(U,G_{j})=0$$by \cite[Proposition 2.2]{shaveisicam}, we have
%$${\rm Tor}_{m+i}^{R}(U,G)\cong{\rm Tor}_{m+i-1}^{R}(U,K_{0})\cong\cdots\cong {\rm Tor}_{i}^{R}(U,N),$$
%as desired.

 Let $G$ be a Gorenstein $(n,d)$-injective $R$-module and ${\rm pd}_{R}(U)=m<\infty$.
Then by hypothesis, the following $(n,d)$-injective resolution of $G$ exists:
$$0\rightarrow N\rightarrow A_{m-1}\rightarrow \cdots \rightarrow A_{0}\rightarrow G\rightarrow 0.$$
So, ${\rm Ext}_{R}^{i}(U, A_j)= 0$ for every $0 \leq j \leq m-1$ and any $i>d$, since $U$  is $n$-presented.  Thus, we deduce that
 ${\rm Ext}_{R}^{i}(U,G)\cong{\rm Ext}_{R}^{m+i}(U, N)=0$.} 

%(4) Setting $G_m=N$ and $K_{j-1}=\ker (G_{j-1}\rightarrow G_{j})$, for every $0\leq j\leq m$, the short exact sequence
%$0\rightarrow K_{j-1}\rightarrow G_{j-1}\rightarrow K_{j}\rightarrow 0$ exist. Thus by (3), the induced exact sequences
%$$0={\rm Ext}_R^{i}(U,G_{j-1})\rightarrow{\rm Ext}_R^{i}(U,K_{j})\rightarrow{\rm Ext}_R^{i+1}(U,K_{j-1})\rightarrow
%{\rm Ext}_R^{i+1}(U,G_{j-1})=0$$
%exists and so ${\rm Ext}_R^{i}(U,K_{j})\cong{\rm Ext}_R^{i+1}(U,K_{j-1})$. Since $K_0=G$, we have
%$${\rm Ext}_R^{m+i}(U,N)\cong{\rm Ext}_R^{m+i-1}(U,K_{m-1})\cong\cdots\cong {\rm Ext}_R^{i}(U,G),$$
%as desired.}
\end{proof}

Next, we study the Gorenstein $(n,d)$-flatness and Gorenstein $(n,d)$-injectivity of modules, in short exact sequences.
\begin{proposition}\label{2.6}
Let $R$ be a left (resp; right) $n$-coherent ring.
\begin{enumerate}
\item [\rm (1)]
Consider the exact sequence $0\rightarrow K\rightarrow B\rightarrow G\rightarrow 0$, where $B$ is $(n,d)$-flat. Then ${\rm G_n^d}$-${\rm fd}(G)\leq{\rm G_n^d}$-${\rm fd}(K)+1$. In particular, if $G$ is Gorenstein $(n,d)$-flat, so is $K$.
\item [\rm (2)]
Consider the exact sequence $0\rightarrow G\rightarrow A\rightarrow N\rightarrow 0$, where $A$ is $(n,d)$-injective. Then ${\rm G_n^d}$-${\rm id}(G)\leq{\rm G_n^d}$-${\rm id}(N)+1$. In particular, if $G$ is Gorenstein $(n,d)$-injective, so is $N$.
\end{enumerate}
\end{proposition}
\begin{proof}
{
(1) We shall show that ${\rm G_n^d}$-${\rm fd}(G)\leq{\rm G_n^d}$-${\rm fd}(K) +1$. In fact, we may assume that
${\rm G_n^d}$-${\rm fd}(K)=m< \infty$. Then, by definition, $K$ admits a Gorenstein $(n,d)$-flat resolution:
$$0\rightarrow B_m\rightarrow B_{m-1}\rightarrow \cdots \rightarrow B_{0}\rightarrow K\rightarrow 0.$$
Assembling this sequence and the short exact sequence $0\rightarrow K\rightarrow B\rightarrow G\rightarrow 0$, the
following commutative diagram is obtained:

\begin{center}
$
\begin{array}{ccccccccccccccccc}
0 &\longrightarrow & B_m & \longrightarrow& \cdots&\longrightarrow &B_1 & \longrightarrow &B_0& \longrightarrow &B& \longrightarrow &G&  \longrightarrow & 0 \\
& & & & &&&&\downarrow & & \uparrow & & \\
 &  &  &  &&&&& K& ={\hspace{-2mm}=}{\hspace{-2mm}={\hspace{-2mm}=}} & K & &  \\
& & & &&&&& \downarrow & & \uparrow & & \\
&& & &&&&& 0 & & 0 & & \\
\end{array}
$
\end{center}
which shows that ${\rm G_n^d}$-${\rm fd}(G)\leq m+1$. The particular case follows from Corollary \ref{2.ty}.

(2) We shall show that ${\rm G_n^d}$-${\rm id}(G)\leq{\rm G_n^d}$-${\rm id}(N) +1$. In fact, we may assume that
${\rm G_n^d}$-${\rm id}(N)=m< \infty$. Then, by definition, $N$ admits a Gorenstein $(n,d)$-injective resolution:
$$0\rightarrow N\rightarrow A_{0}\rightarrow \cdots \rightarrow A_{m-1}\rightarrow A_{m}\rightarrow 0.$$
Assembling this sequence and the short exact sequence $0\rightarrow G\rightarrow A\rightarrow N\rightarrow 0$, the
following commutative diagram is obtained:

\begin{center}
$
\begin{array}{ccccccccccccccccc}
0 &\longrightarrow & G & \longrightarrow & A &\longrightarrow &A_0 & \longrightarrow & A_1& \longrightarrow & \cdots& \longrightarrow &A_m&  \longrightarrow & 0 \\
& & & & \downarrow & & \uparrow & & \\
 &  &  &  & N& ={\hspace{-2mm}=}{\hspace{-2mm}={\hspace{-2mm}=}} & N & &  \\
& & & & \downarrow & & \uparrow & & \\
&& & & 0 & & 0 & & \\
\end{array}
$
\end{center}
which shows that ${\rm G_n^d}$-${\rm id}(G)\leq m+1$. The particular case follows from Corollary \ref{2.h}.}
\end{proof}
\begin{proposition}\label{2.7}
Let $R$ be a left (resp; right)  $n$-coherent ring.
\begin{enumerate}
\item [\rm (1)]
Let $0\rightarrow K\rightarrow G\rightarrow B\rightarrow 0$ be an exact sequence of right $R$-modules. If $K$ is Gorenstein $(n,d)$-flat and $B$ is $(n,d)$-flat, then $G$ is Gorenstein $(n,d)$-flat.
\item [\rm (2)]
Let $0\rightarrow A\rightarrow G\rightarrow N\rightarrow 0$ be an exact sequence of right $R$-modules. If $N$ is Gorenstein $(n,d)$-injective and $A$ is $(n,d)$-injective, then $G$ is Gorenstein $(n,d)$-injective.
\end{enumerate}
\end{proposition}
\begin{proof}
{
(1)
  $K$ is Gorenstein $(n,d)$-flat. So  by Corollary $\ref{2.ty}$, there exists an exact sequence of right $R$-modules
$0\rightarrow K\rightarrow B^{'}\rightarrow L\rightarrow 0$, where $B^{'}$ is $(n,d)$-flat and $L$ is Gorenstein $(n,d)$-flat.
Now, we consider the following diagram:
\begin{center}
$
\begin{array}{ccccccccc}
&&  0 & & 0 & & \\
&&  \downarrow & & \downarrow & & \\
& 0 \longrightarrow & K &\longrightarrow &G & \longrightarrow &B&\longrightarrow &  0 \\
& & \downarrow & & \downarrow& & \parallel & & \\
& 0\longrightarrow & B'& \longrightarrow & N & \longrightarrow&B & \longrightarrow &0 \\
& & \downarrow & & \downarrow & & \\
&&  L & ={\hspace{-1.5mm}=} & L & & \\
& &  \downarrow & & \downarrow & & \\
&&  0 & & 0 & & \\
\end{array}
$
\end{center}
The exactness of the middle horizontal sequence with $B$ and $B^{'}$ are $(n,d)$-flat, implies that $N$ is $(n,d)$-flat. Hence from the middle vertical sequence
and Corollary $\ref{2.ty}$, we deduce that $G$ is Gorenstein $(n,d)$-flat.

(2)  Since $N$ is Gorenstein $(n,d)$-injective, by Corollary $\ref{2.h}$, there exists an exact sequence of of right $R$-modules
$0\rightarrow K\rightarrow A^{'}\rightarrow N\rightarrow 0$, where $A^{'}$ is $(n,d)$-injective and $K$ is Gorenstein $(n,d)$-injective. 
Now, we consider the following diagram:
\begin{center}
$
\begin{array}{ccccccccc}
&& & & 0 & & 0 & & \\
&& & & \downarrow & & \downarrow & & \\
&& & & K & ={\hspace{-1.5mm}=} & K & & \\
&& & & \downarrow & & \downarrow & & \\
0 &\longrightarrow & A & \longrightarrow & D &\longrightarrow &A' & \longrightarrow & 0 \\
& & \parallel& & \downarrow & & \downarrow & & \\
0 & \longrightarrow & A& \longrightarrow & G& \longrightarrow & N & \longrightarrow & 0 \\
& & & & \downarrow & & \downarrow & & \\
&& & & 0 & & 0 & & \\
\end{array}
$
\end{center}
The exactness of the middle horizontal sequence with $A$ and $A^{'}$ are $(n,d)$-injective, implies that $D$ is $(n,d)$-injective. Hence from the middle vertical sequence
and Corollary $\ref{2.h}$, we deduce that $G$ is Gorenstein $(n,d)$-injective.
}
\end{proof}
%%%%%%%%%%%%%%%%%%%%%%%%%%%%%%%%%%%%%%%%%%%%%%%%%%%%%%%%%%%%%%%%%%%%%%%%%%%%%%%%%%%%%%%%%%%%%%%%%%%%%%%%
%%%%%%%%%%%%%%%%%%%%%%%%%%%%%%%%%%%%%%%%%%%%%%%%%%%%%%%%%%%%%%%%%%%%%%%%%%%%%%%%%%%%%%%%%%%%%%%%%%%%%%%%%%%%%%%%%%%%%%%%%%%%%%%%%%%%%%%%%%
\section{ $(n,d)$-Injective  Rings}
A ring $R$ is  right (resp; left)  self-$(n,d)$-injective if $R$ is an $(n,d)$-injective right (resp; left) $R$-module. This section is devoted to $(n,d)$-injective rings over which every $R$-module  is
Gorenstein $(n,d)$-injective.

\begin{proposition}\label{3.1}
Let $R$ be a ring. Then, 
every right $R$-module is Gorenstein $(n,d)$-injective if and only if 
every projective right  $R$-module is $(n,d)$-injective and 
${\rm Ext}_R^{d+1}(U, N)=0$ for any right $R$-module $N$ and any $n$-presented right $R$-module $U$ with ${\rm pd}_{R}(U) <\infty$.
\end{proposition}
\begin{proof}
{
$(\Longrightarrow)$ Let $M$ be a  projective right
$R$-module. Then by (1), $M$ is Gorenstein $(n,d)$-injective. So, the following  $(n,d)$-injective resolution of $M$ exists:
$$\cdots \rightarrow A_1\rightarrow A_0\rightarrow M\rightarrow 0.$$
Since $M$ is projective, $M$ is
$(n,d)$-injective as a direct summand of $A_0$. If $U$ is an $n$-presented right $R$-module with ${\rm pd}_{R}(U) <\infty$, then by Proposition \ref{2.5} and (1),  ${\rm Ext}_R^{d+1}(U, N)=0$ for any right $R$-module $N$. 

$(\Longleftarrow)$ Choose an  injective resolution $0\rightarrow G\rightarrow E^{0}\rightarrow E^{1}\rightarrow\cdots$  and a projective resolution $\cdots \rightarrow F_{1}\rightarrow F_{0}\rightarrow G\rightarrow 0$ of right $R$-module $G$, where every $F_i$ is $(n,d)$-injective by (2). 
 Assembling these resolutions, by Remark \ref{2}, we get the following $(n,d)$-injective resolotion:
$${\mathbf{A}}=\cdots \rightarrow F_{1}\rightarrow F_{0}\rightarrow E^{0}\rightarrow E^{1}\rightarrow\cdots,$$
where $G={\rm ker}(E^0\rightarrow E^1)$, $K^{i}={\rm ker}(E^i\rightarrow E^{i+1})$ and $K_{i}={\rm ker}(F_{i}\rightarrow F_{i-1})$ for any $i\geq 1$.  Let $U$ be a $n$-presented right $R$-module with ${\rm pd}_{R}(U) <\infty$. Then by (2), ${\rm Ext}_R^{d+1}(U, G)={\rm Ext}_R^{d+1}(U, F_i)={\rm Ext}_R^{d+1}(U, E^{i})=0$ for any $i\geq 0$.
So, ${\rm Ext}_R^{d}(U,{\mathbf{A}})$ is exact, and hence $G$ is Gorenstein $(n,d)$-injective.}
\end{proof}

\begin{theorem}\label{3.2}
Let $R$ be a right $n$-coherent ring. Then the following are equivalent:
\begin{enumerate}
\item [\rm (1)]
Every right $R$-module  is Gorenstein $(n,d)$-injective;
\item [\rm (2)]
Every projective right $R$-module is $(n,d)$-injective;
\item [\rm (3)]
$R$ is right self-$(n,d)$-injective.
\end{enumerate}
\end{theorem}
\begin{proof}
{$(1)\Longrightarrow (2)$ and $(2)\Longrightarrow (3)$, is hold by Proposition \ref{3.1}. 

$(3)\Longrightarrow (1)$
Let $G$ be a right $R$-module and $\cdots\rightarrow F_1\rightarrow F_0\rightarrow G\rightarrow 0$ be any
free resolution of $G$. Then by Proposition \ref{3.1},   each $F_i$ is $(n,d)$-injective. Hence Corollary $\ref{2.h}$ completes the proof.}
\end{proof}
\begin{example}\label{3.54} 
{
 Let  $R=k[x^3,x^2,x^{2}y,xy^2,xy,y^2,y^3]$ be  a ring.  We claim that $R$ is not $(n,0)$-injective. Suppose to the contrary, $R$ is $(n,0)$-injective.  $\frac{R}{Rx^2}$ is $n$-presented, since $Rx^2\cong R$ is $n$-presented. Hence $\frac{R}{Rx^2}$ is $n$-presented and ${\rm pd}_{R}(\frac{R}{Rx^2})<\infty$. By Proposition \ref{3.1} and Theorem \ref{3.2}, $\frac{R}{Rx^2}$ is propjective. Threrefore, the exact sequence $0\rightarrow Rx^2\rightarrow R\rightarrow \frac{R}{Rx^2}\rightarrow0$ splits. Thus $Rx^2$ is a direct summand of $R$ and so, $x^2$ is an idempotent, a contradiction.}
\end{example}

\begin{proposition}\label{2.4}
Let $R$ be a ring. Then:
\begin{enumerate}
\item [\rm (1)]
Every Gorenstein $(n,d)$-injective right $R$-module is Gorenstein $(m,d)$-injective for any $m\geq n$.
\item [\rm (2)]
Every Gorenstein $(n,d)$-flat right $R$-module is Gorenstein $(m,d)$-flat for any $m\geq n$.  
\item [\rm (3)]
  If $R$ is a left $n$-coherent, then every Gorenstein $(m,d)$-flat right $R$-module is Gorenstein $(n,d)$-flat for any $m\geq n$.  
\item [\rm (4)]
  If $R$ is a right $n$-coherent, then every Gorenstein $(m,d)$-injective right $R$-module is Gorenstein $(n,d)$-injective for any $m\geq n$.  
\end{enumerate}
\end{proposition}
\begin{proof}
 {(1) Let $G$ be a Gorenstein $(n,d)$-injective right $R$-module. By Remark \ref{2}, every $(n,d)$-injective right $R$-module is $(m,d)$-injective for any $m\geq n$. Thus, there is an exact
sequence
$${\mathbf{I}}= \cdots\longrightarrow I_1\longrightarrow I_{0}\longrightarrow I^0\longrightarrow
I^1\longrightarrow\cdots$$
of $(m,d)$-injective right $R$-modules, where $G={\rm ker}(I^0\rightarrow I^1)$.  By definition and Remark \ref{2}, ${\rm Ext}_{R}^{d}(U,{\mathbf{I}})$ leaves the sequence exact for any $m$-presented right $R$-module $U$ of
finite projective dimension, since every $m$-presented right $R$-module  is $n$-presented and $G$ is Gorenstein $(n,d)$-injective. Hence
$G$ is Gorenstein $(m,d)$-injective.

(4)  Let $G$ be a Gorenstein $(m,d)$-injective right $R$-module. Since $R$  is right $n$-coherent, every $n$-presented right $R$-module  is $m$-presented for any $m\geq n$. Hence, any $(m,d)$-injective right $R$-module is $(n,d)$-injective. Thus, there is an exact
sequence
$${\mathbf{I}}= \cdots\longrightarrow I_1\longrightarrow I_{0}\longrightarrow I^0\longrightarrow
I^1\longrightarrow\cdots$$
of $(n,d)$-injective right $R$-modules, where $G={\rm ker}(I^0\rightarrow I^1)$. So by definition, ${\rm Ext}_{R}^{d}(U,{\mathbf{I}})$ leaves the sequence exact for any $n$-presented right $R$-module $U$ of
finite projective dimension. Hence
$G$ is Gorenstein $(n,d)$-injective. 

(2), (3)  Similar to proof (1), (4).
 }
\end{proof}

A ring $R$ is called right $n$-regular, if every $n$-presented right $R$-module is projective. The following example show that Gorenstein $(m,d)$-injectivity does not imply Gorenstein $(n,d)$-injectivity for any $m\geq n$.
\begin{example}\label{3.94} 
Let $K$ be a field and $E$ be a $k$-vector space with infinite rank. Let $R=K\propto E$. The trivial extension of $K$ by $E$. Then, every right $R$-module is Gorenstein $(2,0)$-injective. We claim that there is right $R$-module $N$ so that $N$ is not Gorenstein $(1,0)$-injective. Suppose to the contrary, every right $R$-module is Gorenstein $(1,0)$-injective. By Proposition \ref{3.1},  we see that ${\rm Ext}_{R}^{1}(U,M)=0$ for any right $R$-module $M$ and any $R$-module $1$-presented $U$. So, every $1$-presented right  $R$-module is projective and it follows that $R$ is $1$-regular or regular, a cotradiction. Since by \cite[Example 3.8]{N.I}, $R$ is not regular.
\end{example}

\begin{theorem}\label{3.83} 
Let $R$ be right and left $n$-coherent ring. Then the following statements are equivalent:
\begin{enumerate}
\item [\rm (1)]
$R$ is two-sided self-$(n,d)$-injective;
\item [\rm (2)]
Every Gorenstein $(n,d)$-flat $R$-module ( right and left) is Gorenstein $(n,d)$-injective;
\item [\rm (3)]
Every Gorenstein flat $R$-module ( right and left) is Gorenstein $(n,d)$-injective;
\item [\rm (4)]
Every flat $R$-module ( right and left) is Gorenstein $(n,d)$-injective;
\item [\rm (5)]
Every Gorenstein projective $R$-module ( right and left) is Gorenstein $(n,d)$-injective;
\item [\rm (6)]
Every projective $R$-module ( right and left) is Gorenstein $(n,d)$-injective;
\item [\rm (7)]
Every Gorenstein injective $R$-module ( right and left) is Gorenstein $(n,d)$-flat;
\item [\rm (8)]
Every injective $R$-module  (right and left) is Gorenstein $(n,d)$-flat;
\item [\rm (9)]
Every Gorenstein $(1,0)$-injective $R$-module ( right and left) is Gorenstein $(n,d)$-flat;
\item [\rm (10)]
Every Gorenstein $(n,0)$-injective $R$-module ( right and left) is Gorenstein $(n,d)$-flat.
\end{enumerate}
\end{theorem}
\begin{proof}
{  $(1)\Longrightarrow (2)$, $(1)\Longrightarrow (3)$, $(1)\Longrightarrow (4)$, $(1)\Longrightarrow (5)$ and $(1)\Longrightarrow (6)$ follow immediately from Theorem $\ref{3.2}$. 

$(3)\Longrightarrow (4)$, $(4)\Longrightarrow (6)$ and $(5)\Longrightarrow (6)$ are trivial.

$(3)\Longrightarrow (1)$ Assume that $G$ is a projective right (resp; left) $R$-module. Then  $G$ is
flat and so, $G$ is Gorenstein
$(n,d)$-injective by (3). So,  similar to proof $(\Longrightarrow)$ of  Proposition $\ref{3.1}$, $G$ is  $(n,d)$-injective. Thus, the assertion follows from Theorem $\ref{3.2}$.

$(6)\Longrightarrow (1)$ This is similar to proof $(3)\Longrightarrow (1)$.

$(1)\Longrightarrow (9)$  By \cite[Proposition 2.21 and Remark 2.22]{N.J}, every $(1,0)$-injective right (resp; left) $R$-module is $(n,d)$-flat. Suppose that $G$ is Gorenstein $(1,0)$-injective. So, the exact sequence
 $${\mathbf{M}}=\cdots \rightarrow M_{1}\rightarrow M_{0}\rightarrow M^{0}\rightarrow M^{1}\rightarrow\cdots,$$
 of $(n,d)$-flat  right (resp; left) $R$-modules exists, where $G={\rm ker}(M^0\rightarrow M^1)$. Let $U$ be a $n$-presented left (resp; right) $R$-module with ${\rm f.d(U)}<\infty$. Then similar to proof Theorem \ref{2.3}(1),  ${\rm Tor}_{d}^{R}(U,{\mathbf{M}})$ is exact. 

$(9)\Longrightarrow (7)$ Let  $M$ be a Gorenstein injective right (resp; left) $R$-module. Then, the exact sequence
 $${\mathbf{M}}=\cdots \rightarrow M_{1}\rightarrow M_{0}\rightarrow M^{0}\rightarrow M^{1}\rightarrow\cdots,$$
 of injective right  (resp; left) $R$-modules exists, where $M={\rm ker}(M^0\rightarrow M^1)$. By Remark \ref{2}, every  injective $R$-module is $(1,0)$-injective. Also, ${\rm Hom}_{R}(U,{\mathbf{M}})$  leaves the sequence exact for any $n$-presented right (resp; left) $R$-module $U$ with ${\rm p.d(U)}<\infty$, and so $M$ is Gorenstein $(1,0)$-injective. 

$(7)\Longrightarrow (8)$
is trivial, since every injective $R$-module is Gorenstein injective.

$(8)\Longrightarrow (1)$
Let $M$ be an injective right $R$-module. Since $M$ is Gorenstein $(n,d)$-flat, we have a long
exact sequence:
 $${\mathbf{M}}=\cdots \rightarrow M_{1}\rightarrow M_{0}\rightarrow M^{0}\rightarrow M^{1}\rightarrow\cdots,$$ 
where any $M_i$ is a $(n,d)$-flat right $R$-module and $M={\rm ker}(M^0\rightarrow M^1)$. Then, the split exact sequence $ 0\rightarrow M\rightarrow M^{0}\rightarrow L\rightarrow 0$ implies that $M$ is $(n,d)$-flat, and hence by \cite[Proposition 2.21 and Remark 2.22]{N.J}, we deduce that $R$ is left self-$(n,d)$-injective. Similarly, $R$ is right self-$(n,d)$-injective.

$(1)\Longrightarrow (10)$
 Suppose that  $G$ is a Gorenstein $(n,0)$-injective right (resp; left) $R$-module. Then by Propsition \ref{2.4}(4), $G$ is  Gorenstein $(1,0)$-injective. Also, by \cite[Proposition 2.21 and Remark 2.22]{N.J}, every $(1,0)$-injective right (resp; left) $R$-module is $(n,d)$-flat. Thus, the exact sequence
 $${\mathbf{N}}=\cdots \rightarrow N_{1}\rightarrow N_{0}\rightarrow N^{0}\rightarrow N^{1}\rightarrow\cdots,$$
 of $(n,d)$-flat right (resp; left) $R$-modules exists, where $G={\rm ker}(N^0\rightarrow N^1)$. Then similar to proof Theorem \ref{2.3}(1),  (10) follows.
 
$(10)\Longrightarrow (9)$ By Proposition \ref{2.4}(1), every Gorenstein $(1,0)$-injective right (resp; left) $R$-module is Gorenstein $(n,0)$-injective. So by (10), (9) is hold. }
\end{proof}
\begin{example}\label{3.54}
{Let $R$ be a commutative $n$-regular ring. Then every $R$-module is Gorenstein $(n,d)$-injective and Gorenstein $(n,d)$-flat, since by \cite[Theorem 3.9]{N.I}, $R$ is self-$(n,0)$-injective and  $n$-coherent. But, there is an  $R_{\infty}$-module $N$ such that $N$ is not Gorenstein $(n,0)$-injective and Gorenstein $(n,0)$-flat. Since if every $R_{\infty}$-module is  Gorenstein $(n,0)$-injective and Gorenstein $(n,0)$-flat, then by Proposition \ref{3.1}, every $n$-presented right $R_{\infty}$-module with finite dimention projective is projective and consequenclay $R_{\infty}$ is $n$-regular and  $n$-coherent, a contradiction. Because $R=K[x^3,x^3y,xy^3,y^3]$ is $n$-regular and if $K$ is a field of characteristic $p>5$, then by \cite[Example 6]{seam3}, $R_{\infty}$ is not $n$-coherent.}
\end{example}
%%%%%%%%%%%%%%%%%%%%%%%%%%%%%%%%%%%%%%%%%%%%%%%%%%%%%%%%%%%%%%%%%%%%%%%%%%%%%%%%%%%%%%%%%%%%%%%%%%%%%%%%
%%%%%%%%%%%%%%%%%%%%%%%%%%%%%%%%%%%%%%%%%%%%%%%%%%%%%%%%%%%%%%%%%%%%%%%%%%%%%%%%%%%%%%%%%%%%%%%%%%%%%%%%%%%%%%%%%%%%%%%%%%%%%%%%%%%%%%%%%%
\section{ Gorenstein $(n, d)$-flat  and
Gorenstein $(n, d)$-injective modules on projective algebras}
\begin{lemma}\label{3.92}
Let $f: R\rightarrow S$ be a surjective ring homomorphism, $S_R$ a
projective right $R$-module  and $_RS$ a
projective left $R$-module. Then:
\begin{enumerate}
\item [\rm (1)]
 If $U$ is an $n$-presented right (resp; left) $R$-module, then $U\otimes_{R}S$ is an $n$-presented right $S$-module.
\item [\rm (2)]
 If $N$ is an $(n,d)$-flat right $R$-module,  then $N\otimes_{R}S$ is an $(n,d)$-flat right $S$-module . 
\item [\rm (3)]
 If $N$ is an $(n,d)$-injective right $R$-module,  then ${\rm Hom}_{R}(S,N)$ is an $(n,d)$-injective right $S$-module . 
\end{enumerate}
\end{lemma}
\begin{proof}
{(1) Suppose  $U$ is $n$-presented right $R$-module. Then,  the exact sequence
 $$0\rightarrow K_{n-1}\rightarrow F_{n-1}\rightarrow\cdots \rightarrow F_{0}\rightarrow F_{1}\rightarrow U\rightarrow 0,$$
of right $R$-modules exists, where $K_{n-1}$ is finitely generated and each $F_{i}$ is finitely generated projective for $0\leq i\leq n-1$.   $_RS$ is  projective, hence the exact sequence
 $$0\rightarrow K_{n-1}\otimes_{R}S\rightarrow F_{n-1}\otimes_{R}S\rightarrow\cdots \rightarrow F_{0}\otimes_{R}S\rightarrow F_{1}\otimes_{R}S\rightarrow U\otimes_{R}S\rightarrow 0,$$
of right $S$-modules exists, where $K_{n-1}\otimes_{R}S$ is finitely generated  and each $F_{i}\otimes_{R}S$ is finitely generated projective for $0\leq i\leq n-1$, since  by \cite[Theorem 2.75]{rotman}, we have 
${\rm Hom}_{R}(F_{i},-)\cong{\rm Hom}_{S}(F_{i}\otimes_{R}S,-)$. So, $U\otimes_{R}S$ is an $n$-presented right $S$-module.

(2) Let $U$ be an $n$-presented left $S$-module. Then by \cite[Lemma 3.10 and Lemma 3.11]{O.H}, $U$ is $n$-presented left $R$-module. By \cite[Corollary 10.61]{rotman}, 
${\rm Tor}_{d+1}^{R}(N,U)\cong{\rm Tor}_{d+1}^{S}(N\otimes_{R}S,U)$, and so $N\otimes_{R}S$ is an $(n,d)$-flat right $S$-module, since  $N$ is  $(n,d)$-flat right  $R$-module. 

(3)  Let $U$ be an $n$-presented right $S$-module. Then $U$ is an $n$-presented right $R$-module. We claim that  ${\rm Ext}_{S}^{k}(U,{\rm Hom}_{R}(S,N))\cong{\rm Ext}_{R}^{k}(U,N)$. We use the induction on $k$. Let $k=0$. Then by \cite[Theorem 2.75]{rotman},${\rm Hom}_{S}(U,{\rm Hom}_{R}(S,N))\cong{\rm Hom}_{R}(U,N)$.  Let $k>0$, then by induction hypothesis and the exact sequence $0\rightarrow K\rightarrow F_{0}\rightarrow U \rightarrow 0$, where $F_0$ is a finitely generated free,  we deduce that
 ${\rm Ext}_{S}^{k}(U,{\rm Hom}_{R}(S,N))\cong{\rm Ext}_{R}^{k}(U,N)$, and so, (3) is follows.
}
\end{proof}
\begin{proposition}\label{3.95}
Let $f: R\rightarrow S$ be a surjective ring homomorphism, $S_R$ a
projective right $R$-module  and $_RS$ a
projective left $R$-module. Then:
\begin{enumerate}
\item [\rm (1)]
If $U$ is a left $S$-module with ${\rm fd}_S(U)<\infty$, then ${\rm fd}_R(U)<\infty$.
\item [\rm (2)]
If $U$ is a left $R$-module with ${\rm fd}_R(U)<\infty$, then ${\rm fd}_S(U\otimes_{R}S)<\infty$.
\item [\rm (3)]
If $U$ is a right $S$-module with ${\rm pd}_S(U)<\infty$, then ${\rm pd}_R(U)<\infty$.
\item [\rm (4)]
If $U$ is a right $R$-module with ${\rm pd}_R(U)<\infty$, then ${\rm pd}_S(U\otimes_{R}S)<\infty$.
\end{enumerate}
\end{proposition}
\begin{proof}
{(1) By \cite[Lemma 3.10]{O.H} and  \cite[Corollary 10.61]{rotman}, 
$${\rm Tor}_{k}^{R}(-,U)\cong{\rm Tor}_{k}^{R}(-,U\otimes_{R}S)\cong{\rm Tor}_{k}^{S}(-,U)$$. So, ${\rm fd}_R(U)<\infty$.

(2) By \cite[Corollary 10.61]{rotman},  ${\rm Tor}_{n}^{R}(U,-)\cong{\rm Tor}_{n}^{S}(U\otimes_{R}S,-)$. So, (2) holds.

(3) Similar to proof (3) of lemma \ref{3.92},  ${\rm Ext}_{S}^{n}(U,{\rm Hom}_{R}(S,-))\cong{\rm Ext}_{R}^{n}(U,-).$ So, (2) holds.

(4) By \cite[Theorem 2.76]{rotman}, 
${\rm Ext}_{R}^{n}(U.-)\cong{\rm Ext}_{S}^{n}(U\otimes_{R}S.-)$. So, ${\rm pd}_S(U\otimes_{R}S)<\infty$.

}
\end{proof}
\begin{theorem}\label{3.90}
Let $f: R\rightarrow S$ be a surjective ring homomorphism, $S_R$ a
projective right $R$-module,  $_RS$ a
projective left $R$-module and $G$ a left $S$-module. Then, the following statements are equivalent:
\begin{enumerate}
\item [\rm (1)]
 $G_{R}$ is Gorenstein $(n,d)$-flat;
\item [\rm (2)] $G\otimes_{R}S$ is Gorenstein $(n,d)$-flat;
\item [\rm (3)]
$G_{S}$ is Gorenstein $(n,d)$-flat.
\end{enumerate}
\end{theorem}
\begin{proof}
{$(1)\Longrightarrow (2)$  Suppose that  $N$ is an $(n,d)$-flat right $R$-module. Then by Lemma \ref{3.92},  $N\otimes_{R}S$ is an $(n,d)$-flat right $S$-module. But, $_RS_{S}$ is a flat. So if  $G$ is an Gorenstein $(n,d)$-flat right $R$-module, then the exact sequence
 $${\mathbf{F}}=\cdots \rightarrow F_{1}\rightarrow F_{0}\rightarrow F^{0}\rightarrow F^{1}\rightarrow\cdots,$$
 of $(n,d)$-flat right $R$-modules,  where $G={\rm ker}(F^0\rightarrow F^1)$,  induces the following exact sequence of $(n,d)$-flat  right $S$-modules :
 $${\mathbf{F}}\otimes_{R}S=\cdots \rightarrow F_{1}\otimes_{R}S\rightarrow F_{0}\otimes_{R}S\rightarrow F^{0}\otimes_{R}S\rightarrow F^{1}\otimes_{R}S\rightarrow\cdots,$$
 where $G\otimes_{R}S={\rm ker}(F^0\otimes_{R}S\rightarrow F^1\otimes_{R}S)$. Let $U$ be an $n$-presented left $S$-module with ${\rm fd}_{S}(U)<\infty$. We show that ${\rm Tor}_{d}^{S}(U,{\mathbf{F}}\otimes_{R}S) $ is exact.  By  Proposition \ref{3.95}, ${\rm fd}_{R}(U)<\infty$.  Also, $U$ is an $n$-presented left $R$-module. Thus by hypothesis, ${\rm Tor}_{d}^{R}(U,{\mathbf{F}})$ is exact. Therefore ${\rm Tor}_{d}^{S}(U,{\mathbf{F}}\otimes_{R}S) $ is exact, since by \cite[Corollary 10.61]{rotman}, ${\rm Tor}_{d}^{R}(U,{\mathbf{F}})\cong{\rm Tor}_{d}^{R}(U\otimes_{R}S,{\mathbf{F}})\cong {\rm Tor}_{d}^{S}(U,{\mathbf{F}}\otimes_{R}S) $. Hence, $G\otimes_{R}S$ is Gorenstein $(n,d)$-flat.

$(2)\Longrightarrow (3)$ By \cite[Lemma 3.10]{O.H}, $G\otimes_{R}S\cong G_{S}$, and hence  $G_{S}$ is  Gorenstein $(n,d)$-flat.

$(3)\Longrightarrow (1)$ Suppose that  $N$ is an $(n,d)$-flat right $S$-module. Then by \cite[Lemma 3.12]{O.H},  $N$ is an $(n,d)$-flat right $R$-module. So, if   $G$ be a Gorenstein $(n,d)$-flat right $S$-module, then the exact sequence
 $${\mathbf{F}}=\cdots \rightarrow F_{1}\rightarrow F_{0}\rightarrow F^{0}\rightarrow F^{1}\rightarrow\cdots,$$
 of $(n,d)$-flat right $R$-modules exists, where $G={\rm ker}(F^0\rightarrow F^1)$. Let $U$ be an $n$-presented left $R$-module with ${\rm fd}_{R}(U)<\infty$. It is sufficient to prove that ${\rm Tor}_{d}^{R}(U,{\mathbf{F}})$  is exact.  $U\otimes_{R}S$ is an $n$-presented left $S$-module by Lemma \ref{3.92}, and by Proposition \ref{3.95}, ${\rm fd}_{S}(U\otimes_{R}S)<\infty$.
 So by hypothesis, ${\rm Tor}_{d}^{S}(U\otimes_{R}S,{\mathbf{F}})$ is exact.  On the other hand, by \cite[Corollary 10.61]{rotman}, we have ${\rm Tor}_{d}^{S}(U\otimes_{R}S,{\mathbf{F}})\cong {\rm Tor}_{d}^{R}(U,{\mathbf{F}}\otimes_{S}S)\cong{\rm Tor}_{d}^{R}(U,{\mathbf{F}})$.  Therefore, ${\rm Tor}_{d}^{R}(U,{\mathbf{F}})$  is exact, and hence  $G_{R}$ is Gorenstein $(n,d)$-flat.

}
\end{proof}
\begin{theorem}\label{3.96}
Let $f: R\rightarrow S$ be a surjective ring homomorphism, $S_R$ a
projective right $R$-module,  $_RS$ a
projective left $R$-module and $G$ a left $S$-module. Then,  the following statements are equivalent:
\begin{enumerate}
\item [\rm (1)]
 $G_{R}$ is  Gorenstein $(n,d)$-injective;
\item [\rm (2)]
 ${\rm Hom}_{R}(S,G)$ is Gorenstein $(n,d)$-injective;
\item [\rm (3)]
 $G_{S}$ is Gorenstein $(n,d)$-injective.
\end{enumerate}
\end{theorem}
\begin{proof}
{$(1)\Longrightarrow (2)$
Suppose that  $N$ is an $(n,d)$-injective right $R$-module. Then by Lemma \ref{3.92},  ${\rm Hom}_{S}(S,N)$ is an $(n,d)$-injective right $S$-module. But, $S$ is a projective $R$-module. So if $G$ is Gorenstein $(n,d)$-injective right $R$-module , then the exact sequence
 $${\mathbf{A}}=\cdots \rightarrow A_{1}\rightarrow A_{0}\rightarrow A^{0}\rightarrow A^{1}\rightarrow\cdots,$$
 of $(n,d)$-injective right $R$-modules,  where $G={\rm ker}(A^0\rightarrow A^1)$,  induces the following exact sequence of $(n,d)$-injective right $S$-modules :
 $${\rm Hom}_{R}(S,{\mathbf{A}})=\cdots \rightarrow {\rm Hom}_{R}(S,A_{1})\rightarrow {\rm Hom}_{R}(S,A_{0})\rightarrow {\rm Hom}_{R}(S,A^{0})\rightarrow \cdots,$$
 where ${\rm Hom}_{R}(S,G)={\rm ker}({\rm Hom}_{R}(S,A^0)\rightarrow {\rm Hom}_{R}(S,A^1))$. Let $U$ be an $n$-presented right $S$-module with ${\rm pd}_{S}(U)<\infty$. Then by \cite[Lemmas 3.10 and 3.11]{O.H}, $U\cong U\otimes_{R}S$ is an $n$-presented right $R$-module, and by  Proposition \ref{3.95}, ${\rm pd}_{R}(U)<\infty$.   Also, similar to proof (3) of Lemma \ref{3.92}, we have that $${\rm Ext}_{S}^{d}(U,{\rm Hom}_{R}(S,{\mathbf{A}}))\cong{\rm Ext}_{R}^{d}(U,{\mathbf{A}}).$$
By hypothesis, ${\rm Ext}_{R}^{d}(U,{\mathbf{A}})$ is exact. Therefore ${\rm Ext}_{S}^{d}(U,{\rm Hom}_{R}(S,{\mathbf{A}}))$ is exact, and hence ${\rm Hom}_{R}(S,G)$ is Gorenstein $(n,d)$-injective.

$(2)\Longrightarrow (3)$ By \cite[Proposition 8.33]{rotman}, 
${\rm Hom}_{R}(S,G)={\rm Hom}_{S}(S,G)\cong G$, since $f$ is surjective. So, we deduce that  $G_{S}$ is Gorenstein $(n,d)$-injective.

$(3)\Longrightarrow (1)$
Suppose that  $N$ is an $(n,d)$-injective right $S$-module. Then by \cite[Lemma 3.12]{O.H},  $N$ is an $(n,d)$-injective right $R$-module. So, if   $G$ be a Gorenstein $(n,d)$-injective right $S$-module, then the exact sequence
 $${\mathbf{A}}=\cdots \rightarrow A_{1}\rightarrow A_{0}\rightarrow A^{0}\rightarrow A^{1}\rightarrow\cdots,$$
 of $(n,d)$-injective right $R$-modules exists, where $G={\rm ker}(A^0\rightarrow A^1)$. Let $U$ be an $n$-presented right $R$-module with ${\rm pd}_{R}(U)<\infty$. It is sufficient to prove that ${\rm Ext}_{R}^{d}(U,{\mathbf{A}})$ is exact.  $U\otimes_{R}S$ is an $n$-presented right $S$-module by Lemma \ref{3.92}. Also by Proposition \ref{3.95}, ${\rm pd}_{S}(U\otimes_{R}S)<\infty$. So by hypothesis, ${\rm Ext}_{S}^{d}(U\otimes_{R}S, {\mathbf{A}})$ is exact. On the other hand, similar to proof (4) of Proposition \ref{3.95}, we have ${\rm Ext}_{R}^{d}(U,{\mathbf{A}})\cong{\rm Ext}_{S}^{d}(U\otimes_{R}S, {\mathbf{A}})$. Therefore, ${\rm Ext}_{R}^{d}(U,{\mathbf{A}})$  is exact, and hence  $G_{R}$ is a Gorenstein $(n,d)$-injective.
}
\end{proof}

we have the following corollary.
\begin{corollary}\label{3.98}
Let $f: R\rightarrow S$ be a surjective ring homomorphism, $S$ a
projective $R$-module and $M$ a $S$-module. Then :
\begin{enumerate}
\item [\rm (1)]
$G_{d}^{n}.fd(M_R) = G_{d}^{n}.fd(M_S) = G_{d}^{n}.fd(S\otimes_{R}M).$
\item [\rm (2)]
 $G_{d}^{n}.id(M_R) = G_{d}^{n}.id(M_S) = G_{d}^{n}.id({\rm Hom}_{R}(S,M)).$
\end{enumerate}
\end{corollary}
     \  \  \   In the following propositon, Gorenstein $(n,d)$-injectivity of modules over almost excellent extensions is studied.
 \begin{proposition}\label{3.nv}
Let  $S\geq R$ be an almost excellent extension. Then
\begin{enumerate}
\item [\rm (1)]
 $R$ is right Gorenstein $(n,d)$-injective if and only if right $S$ is Gorenstein $(n,d)$-injective.
\item [\rm (2)]
 $R_{\infty}$ is right Gorenstein $(n,d)$-injective if and only if  $S_{\infty}$ is right Gorenstein $(n,d)$-injective.
\end{enumerate}
\end{proposition}
\begin{proof}
{(1) ($\Longrightarrow$): Let $R$ be a right Gorenstein $(n,d)$-injective, then the exact sequence
 $${\mathbf{A}}=\cdots \rightarrow A_{1}\rightarrow A_{0}\rightarrow A^{0}\rightarrow A^{1}\rightarrow\cdots,$$
 of $(n,d)$-injective right $R$-modules exists, where $R={\rm ker}(A^0\rightarrow A^1)$. Since the  exact sequence $0\rightarrow K\rightarrow A_{0}\rightarrow R\rightarrow 0 $ is split, $R$ is right $(n,d)$-injective. Thus by \cite[Proposition 5.1]{X.J}, $S$ is right $(n,d)$-injective and hence by Remark \ref{2}, $S$ is right Gorenstein $(n,d)$-injective.

 ($\Longleftarrow$):  Similar to proof ($\Longrightarrow$).

(2) If $S\geq R$ is an almost excellent extension, then $S_{\infty}\geq R_{\infty}$ is an almost excellent extension. Hence by (1), (2) is hold.
}
\end{proof}
From Theorem \ref{3.83}, Proposition \ref{3.nv} and \cite[Theorem 5.3]{X.J} we immediately have the following corollary.
 \begin{corollary}\label{3.np}
Let  $S\geq R$ be an almost excellent extension and $R$ be $n$-coherent. Then the following statements are equivalent:
\begin{enumerate}
\item [\rm (1)]
 $R_R$ is Gorenstein $(n,d)$-injective;
\item [\rm (2)]
Every Gorenstein $(n,d)$-flat right $R$-module is Gorenstein $(n,d)$-injective right $R$-module;
\item [\rm (3)]
Every Gorenstein $(n,d)$-flat right $S$-module is Gorenstein $(n,d)$-injective right $S$-module.
\end{enumerate}
\end{corollary}
%{\bf Acknowledgements.} We thank the referee(s) for helpful suggestions on an earlier
%version of this paper.

%\section*{Acknowledgments}
%The author wish to thank $\cdots$

%%%%%%%%%%%%%%%%%%%%%%%%%%%%%%%%%%%%%%%%%%%%%%%%%%%%%%%%%%
\providecommand{\bysame}{\leavevmode\hbox
to3em{\hrulefill}\thinspace}

%%%%%%%%%%%%%%%%%%%%%%%%%%%%%%%%%%%%%%%%%%%%%%%%%%%%%%%%%%%%%%%%%%%%%%%%%

\end{document}